# ONE-DEPENDENT TRIGONOMETRIC DETERMINANTAL PROCESSES ARE TWO-BLOCK-FACTORS[1]

BY ERIK I. BROMAN

*Chalmers University of Technology*

Given a trigonometric polynomial $f:[0,1] \to [0,1]$ of degree $m$, one can define a corresponding stationary process $\{X_i\}_{i \in \mathbb{Z}}$ via determinants of the Toeplitz matrix for $f$. We show that for $m = 1$ this process, which is trivially one-dependent, is a two-block-factor.

**1. Introduction.** We will start by defining a family of probability measures $\mathbf{P}^f$ on the Borel sets of $\{0,1\}^{\mathbb{Z}}$ where $f:[0,1] \to [0,1]$ is a Lebesgue-measurable function (see [9]). For such an $f$, define the probability of the cylinder sets by

$$\mathbf{P}^f[\eta(e_1) = \cdots = \eta(e_k) = 1]$$
$$:= \mathbf{P}^f[\{\eta \in \{0,1\}^{\mathbb{Z}} : \eta(e_1) = \cdots = \eta(e_k) = 1\}]$$
$$:= \det[\hat{f}(e_j - e_i)]_{1 \leq i,j \leq k},$$

where $e_1, \ldots, e_k$ are distinct elements in $\mathbb{Z}$ and $k \geq 1$. Here $\hat{f}$ denotes the Fourier coefficients of $f$, defined by

$$\hat{f}(k) := \int_0^1 f(x) e^{-i2\pi k x} \, dx.$$

In [9] it is proven that $\mathbf{P}^f$ is indeed a probability measure. In fact, they showed this for the more general case of $f:\mathbb{T}^d \to [0,1]$, where $\mathbb{T}^d := \mathbb{R}^d/\mathbb{Z}^d$; in this case, the resulting process is indexed by $\mathbb{Z}^d$. This result rests very strongly on the results in [8]. Except for the two definitions below, $\{X_i\}_{i \in \mathbb{Z}}$ will always denote a process distributed according to some measure $\mathbf{P}^f$. Throughout this paper, equality in distribution will be denoted by $\stackrel{\mathcal{D}}{=}$. Let

Received September 2003; revised September 2003.

[1]Supported in part by the Swedish Natural Science Research Council.

the function $f:[0,1] \to [0,1]$ be of the form

$$f(x) = \sum_{k=-m}^{m} a_k e^{-i2\pi kx}.$$

It is then easily checked that the process $\{X_i\}_{i \in \mathbb{Z}}$ corresponding to the probability measure $\mathbf{P}^f$ is $m$-dependent according to the definition below.

DEFINITION 1.1. A process $\{X_i\}_{i \in \mathbb{Z}}$ is called $m$-dependent if $\{X_i\}_{i<k}$ is independent of $\{X_i\}_{i \geq k+m}$ for all integers $k$.

We will also need the definition of an $m$-block-factor.

DEFINITION 1.2. The process $\{X_i\}_{i \in \mathbb{Z}}$ is called an $m$-block-factor if there exists a function $h$ of $m$ variables and an i.i.d. process $\{Y_i\}_{i \in \mathbb{Z}}$ such that $\{X_i\}_{i \in \mathbb{Z}} \stackrel{\mathcal{D}}{=} \{h(Y_i, \ldots, Y_{i+m-1})\}_{i \in \mathbb{Z}}$.

We will, as usual, not distinguish between the process $\{X_i\}_{i \in \mathbb{Z}}$ and the corresponding probability measure $\mathbf{P}^f$.

Observe that an $(m+1)$-block-factor is trivially $m$-dependent. For some time, it was an open question whether all $m$-dependent processes were in fact $(m+1)$-block-factors (see [[4]–[7]]). However, in [2] the authors constructed a family of one-dependent processes which are not two-block-factors, and in [3] the authors constructed a one-dependent process which is not a $k$-block-factor for any $k$. In [1] the authors constructed a one-dependent stationary Markov process with five states which is not a two-block-factor; they also proved that this result is sharp in the sense that every one-dependent stationary Markov process with not more than four states is a two-block-factor. In view of the above, it is a natural question to ask whether a certain $m$-dependent process is an $(m+1)$-block-factor or not.

$\mathbf{P}^f$ as defined above is an $m$-dependent "trigonometric determinantal probability measure." These probability measures are special cases of general determinantal probability measures; see [8] or [10] for definitions and results. Determinantal processes arise in numerous contexts, for example, mathematical physics, random matrix theory and representation theory, to name a few. For a survey, see [10], for further results, see [8] and for results concerning the discrete stationary case, see [9]. In [9], they ask whether $\mathbf{P}^f$ above is an $(m+1)$-block-factor. In that paper they say that if one can find sufficiently explicit block factors for all trigonometric polynomials, then one can find explicit factors of i.i.d. processes giving $\mathbf{P}^f$, where $f$ is any function such that $f:\mathbb{T} \to [0,1]$. This in turn would enable one to use more standard probabilistic techniques when studying such a $\mathbf{P}^f$. We answer their question positively for $m=1$ in Theorem 1.3, constructing an explicit two-block-factor.

THEOREM 1.3. *If $f:[0,1] \to [0,1]$ is given by*

$$f(x) = b + ae^{-i2\pi x} + ce^{i2\pi x},$$

*then the corresponding process $\{X_i\}_{i \in \mathbb{Z}}$ is a two-block-factor.*

**2. Proof of Theorem 1.3.** With $f$ as in the statement of the theorem, it follows that $\bar{a} = c$, $b \geq 0$ and hence if $a = a_1 + ia_2$,

(1) $\quad f(x) = b + 2a_1\cos(2\pi x) + 2a_2\sin(2\pi x) = b + 2|a|\cos(2\pi x - \phi),$

for some suitable choice of $\phi$. Let, as usual, $\mathbf{P}^f$ be the corresponding probability measure, and write

$$D_k := \det[\hat{f}(j-i)]_{1 \leq i,j \leq k+1},$$

where $k \geq 0$.

Note that the process $\{X_i\}_{i \in \mathbb{Z}}$ distributed according to $\mathbf{P}^f$ is obviously stationary. Since $\mathbf{P}^f$ is one-dependent, it is easily seen that it is uniquely determined among the one-dependent processes by the values of

$$\mathbf{P}^f[\eta(i) = \cdots = \eta(i+k) = 1] = \mathbf{P}^f[\eta(1) = \cdots = \eta(1+k) = 1]$$

as $k$ varies over the nonnegative integers.

We have that for $k \geq 2$,

(2)
$$D_k = \det[\hat{f}(j-i)]_{1 \leq i,j \leq k+1} = \begin{vmatrix} b & a & 0 & 0 & 0 & \cdots \\ \bar{a} & b & a & 0 & 0 & \cdots \\ 0 & \bar{a} & b & a & 0 & \cdots \\ 0 & 0 & \bar{a} & b & a & \cdots \\ \vdots & \vdots & \vdots & \ddots & \ddots & \ddots \end{vmatrix}$$

$$= b \begin{vmatrix} b & a & 0 & 0 & \cdots \\ \bar{a} & b & a & 0 & \cdots \\ 0 & \bar{a} & b & a & \cdots \\ 0 & 0 & \bar{a} & b & \cdots \\ \vdots & \vdots & \vdots & \ddots & \ddots \end{vmatrix} - a \begin{vmatrix} \bar{a} & a & 0 & 0 & \cdots \\ 0 & b & a & 0 & \cdots \\ 0 & \bar{a} & b & a & \cdots \\ 0 & 0 & \bar{a} & b & \cdots \\ \vdots & \vdots & \vdots & \ddots & \ddots \end{vmatrix}$$

$$= bD_{k-1} - |a|^2 D_{k-2},$$

where the determinant on the left-hand side of the third equality has size $(k+1) \times (k+1)$, and the two on the right-hand side have size $k \times k$. Furthermore,

(3) $\qquad\qquad\qquad D_0 = |b| = b,$

(4) $\qquad\qquad\qquad D_1 = \begin{vmatrix} b & a \\ \bar{a} & b \end{vmatrix} = b^2 - |a|^2.$

The characteristic equation corresponding to (2) is

(5) $$r^2 - br + |a|^2 = 0,$$

which has two roots

(6) $$r_1 = \frac{b}{2} + \sqrt{\frac{b^2}{4} - |a|^2}$$

and

(7) $$r_2 = \frac{b}{2} - \sqrt{\frac{b^2}{4} - |a|^2}.$$

*Case* 1. Assume that $r_1 = r_2 = r$ so that $r = \frac{b}{2}$ and

$$\frac{b^2}{4} = |a|^2,$$

and so (since $b, |a| \geq 0$)

$$b = 2|a|.$$

We have by (1) that

$$\max_{x \in [0,1]} f(x) = \max_{x \in [0,1]} (b + 2|a|\cos(2\pi x - \phi)) = b + 2|a| = 2b,$$

and since $f:[0,1] \to [0,1]$, we get $b \leq \frac{1}{2}$ and so $|a| \leq \frac{1}{4}$.

With $r_1 = r_2 = r$, it follows from the basic theory of difference equations that the solution to (2) is

$$D_k = (C_1 k + C_2) r^k \qquad \forall k \geq 0,$$

for some constants $C_1, C_2$ yet to be determined. Using (3) and (4), we get that $C_2 = D_0 = b = 2r$, and using this, we get $(C_1 + 2r)r = D_1 = b^2 - |a|^2 = b^2 - b^2/4 = 3r^2$. Hence $C_1 = r$ and so

(8) $$D_k = (kr + 2r) r^k \qquad \forall k \geq 0.$$

We will now construct a two-block-factor which we will show to be distributed according to $\mathbf{P}^f$. Let $\{Y_i\}_{i \in \mathbb{Z}}$ be i.i.d. uniform on $[0,1]$. Define $h:[0,1] \times [0,1] \to [0,1]$ by $h = I_A$, where

$$A = [0, \tfrac{1}{4}] \times [0, r] \cup [0, \tfrac{1}{4}] \times [\tfrac{1}{2}, \tfrac{1}{2} + r] \cup [0, \tfrac{1}{4}] \times [\tfrac{3}{4}, \tfrac{3}{4} + r]$$
$$\cup [\tfrac{1}{4}, \tfrac{1}{2}] \times [\tfrac{1}{4}, \tfrac{1}{4} + r] \cup [\tfrac{1}{4}, \tfrac{1}{2}] \times [\tfrac{1}{2}, \tfrac{1}{2} + r] \cup [\tfrac{1}{4}, \tfrac{1}{2}] \times [\tfrac{3}{4}, \tfrac{3}{4} + r]$$
$$\cup [\tfrac{1}{2}, \tfrac{3}{4}] \times [\tfrac{1}{2}, \tfrac{1}{2} + r] \cup [\tfrac{3}{4}, 1] \times [\tfrac{3}{4}, \tfrac{3}{4} + r].$$

$A$ is depicted as the gray area of Figure 1. Observe that $r = |a| \leq \frac{1}{4}$.

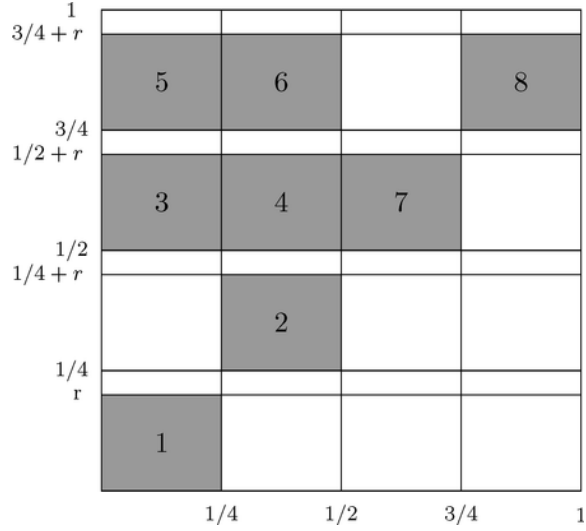

Fig. 1. *This figure shows A (the shaded area).*

We will show that

$$\mathbf{P}[h(Y_i, Y_{i+1}) = \cdots = h(Y_{i+k}, Y_{i+k+1}) = 1] = D_k \qquad \forall\, k \geq 0.$$

Since $\{h(Y_i, Y_{i+1})\}_{i \in \mathbb{Z}}$ is one-dependent, this gives us $\{h(Y_i, Y_{i+1})\}_{i \in \mathbb{Z}} \stackrel{\mathcal{D}}{=} \mathbf{P}^f$ as desired. We first observe that the size of the shaded area of Figure 1 is $8\frac{1}{4}r = 2r = b$, so that $\mathbf{P}[h(Y_i, Y_{i+1}) = 1] = D_0$.

If $h(Y_i, Y_{i+1}) = \cdots = h(Y_{i+k}, Y_{i+k+1}) = 1$, then $(Y_{i+l}, Y_{i+l+1})$ must be in one of the boxes marked 1 through 8 of Figure 1 $\forall\, l \in \{0, \ldots, k\}$. If $(Y_i, Y_{i+1})$ is in the box marked 1, then $Y_{i+1} \in [0, r]$ and so $(Y_{i+1}, Y_{i+2})$ must be in one of the boxes marked 1, 3 or 5 because otherwise $(Y_{i+1}, Y_{i+2}) \notin A$. Similar "rules" apply if $(Y_i, Y_{i+1})$ is in one of the other seven boxes. We see that for any $\omega$ such that $h(Y_i(\omega), Y_{i+1}(\omega)) = \cdots = h(Y_{i+k}(\omega), Y_{i+k+1}(\omega)) = 1$, there is a natural sequence, $j_0 j_1 \cdots j_k(\omega) \in \{1, \ldots, 8\}^{k+1}$ associated to it, where the value of $j_l$ indicates that $(Y_{i+l}(\omega), Y_{i+l+1}(\omega))$ is in the box marked with that value. In any such sequence, the number 1 can only be followed by 1, 3 or 5, as described above, while the number 2 can only be followed by 2, 4 or 6. Additionally, any one of the numbers 3, 4 or 7 must be followed by a 7, while any one of 5, 6 or 8 must be followed by an 8.

We claim that the number of sequences $j_0 j_1 \cdots j_k$ described above is $(4k + 8)$. To see this, observe that every such sequence with $j_k \notin \{1, 2\}$ can be extended into a sequence $j_0 j_1 \cdots j_{k+1}$ in only one way, while if $j_k \in \{1, 2\}$, it can be extended in three ways. Observe also that there are only two sequences $j_0 j_1 \cdots j_k$ ending in 1 or 2.

The set of $\omega$ giving a specific sequence $j_0 j_1 \cdots j_k \in \{1, \ldots, 8\}^{k+1}$ has probability $\frac{1}{4} r^{k+1}$ since $Y_i$ must be in an interval of length $\frac{1}{4}$, while $Y_{i+1}, \ldots, Y_{i+k+1}$ all must be within intervals of length $r$. Hence the total probability of having $h(Y_i, Y_{i+1}) = \cdots = h(Y_{i+k}, Y_{i+k+1}) = 1$ is $(4k+8)\frac{1}{4}r^{k+1} = (kr+2r)r^k$. Comparing with (8), we see that

$$\mathbf{P}[h(Y_i, Y_{i+1}) = \cdots = h(Y_{i+k}, Y_{i+k+1}) = 1] = D_k$$

$\forall k \geq 0$, and we conclude that $\{h(Y_i, Y_{i+1})\}_{i \in \mathbb{Z}} \stackrel{\mathcal{D}}{=} \mathbf{P}^f$ and so this case is proved.

*Case* 2. It remains to consider $r_1 \neq r_2$. According to (6) and (7),

$$r_1 + r_2 = b$$

and

$$r_1 r_2 = |a|^2.$$

In this case the solution to (2) is, again, from basic difference equation theory,

$$D_k = C_1 r_1^k + C_2 r_2^k \qquad \forall k \geq 0,$$

for some constants $C_1, C_2$ yet to be determined. Using this with (3), we get

$$C_1 + C_2 = D_0 = r_1 + r_2,$$

and using (4), we get

$$C_1 r_1 + C_2 r_2 = D_1 = b^2 - |a|^2 = (r_1 + r_2)^2 - r_1 r_2 = r_1^2 + r_1 r_2 + r_2^2.$$

A straightforward calculation yields

$$C_1 = \frac{r_1^2}{r_1 - r_2}$$

and

$$C_2 = -\frac{r_2^2}{r_1 - r_2},$$

and therefore for $k \geq 1$,

$$D_k = \frac{r_1^{k+2} - r_2^{k+2}}{r_1 - r_2} = \frac{r_1^{k+2} - r_1^{k+1} r_2 + r_2(r_1^{k+1} - r_2^{k+1})}{r_1 - r_2} = r_1^{k+1} + r_2 D_{k-1}.$$
(9)

Assume that $b \leq \frac{1}{2}$ so that $2(r_1 + r_2) \leq 1$. We will now construct a two-block-factor which we will show to be distributed according to $\mathbf{P}^f$. Let $\{Y_i\}_{i \in \mathbb{Z}}$ be i.i.d. uniform on $[0,1]$ and again take $h\colon [0,1]\times[0,1]\to[0,1]$ to be the function $h = I_A$, where $A$ is now

$$\begin{aligned}A = {} & [0, Cr_1]\times[0, r_1] \cup [0, Cr_1]\times[2Cr_1, 2Cr_1+r_2]\\ & \cup [0, Cr_1]\times[2Cr_1+Cr_2, 2Cr_1+Cr_2+r_2]\\ & \cup [Cr_1, 2Cr_1]\times[Cr_1, Cr_1+r_1]\\ & \cup [Cr_1, 2Cr_1]\times[2Cr_1, 2Cr_1+r_2]\\ & \cup [Cr_1, 2Cr_1]\times[2Cr_1+Cr_2, 2Cr_1+Cr_2+r_2]\\ & \cup [2Cr_1, 2Cr_1+Cr_2]\times[2Cr_1, 2Cr_1+r_2]\\ & \cup [2Cr_1+Cr_2, 1]\times[2Cr_1+Cr_2, 2Cr_1+Cr_2+r_2]\end{aligned}$$

and $C = \frac{1}{2(r_1+r_2)} \geq 1$. $A$ is the shaded area of Figure 2.

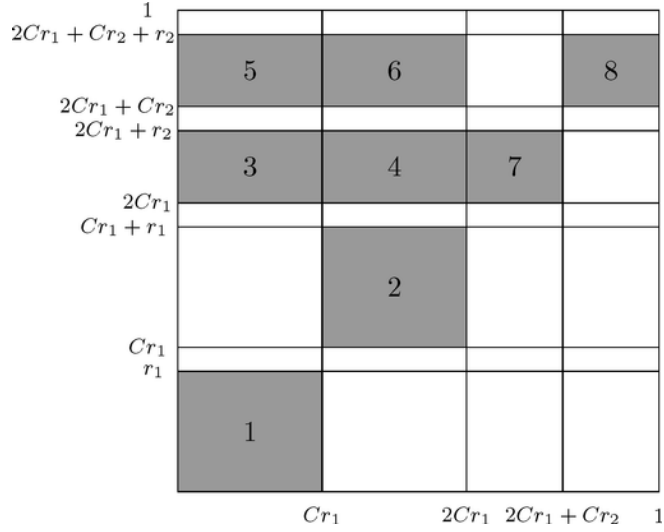

Fig. 2. *This figure shows A (the shaded area).*

Again we will show that

$$\mathbf{P}[h(Y_i, Y_{i+1}) = \cdots = h(Y_{i+k}, Y_{i+k+1}) = 1] = D_k \qquad \forall\, k \geq 0.$$

Since again $\{h(Y_i, Y_{i+1})\}_{i\in\mathbb{Z}}$ is one-dependent, this gives us $\{h(Y_i, Y_{i+1})\}_{i\in\mathbb{Z}} \stackrel{\mathcal{D}}{=} \mathbf{P}^f$. We observe that the size of the shaded area of Figure 2 equals

$$2Cr_1r_1 + 4Cr_1r_2 + 2Cr_2r_2 = 2C(r_1+r_2)^2 = r_1+r_2$$

by our choice of $C$, and so $\mathbf{P}[h(Y_i, Y_{i+1}) = 1] = D_0$.

For any $\omega$ such that $h(Y_i(\omega), Y_{i+1}(\omega)) = \cdots = h(Y_{i+k}(\omega), Y_{i+k+1}(\omega)) = 1$, there is a natural sequence $j_0 j_1 \cdots j_k(\omega) \in \{1,\ldots,8\}^{k+1}$ associated to it, as before. Let $\{\omega : j_0 j_1 \cdots j_k(\omega)\}$ denote the set of $\omega$ giving a specific sequence $j_0 j_1 \cdots j_k$, and for convenience we will write $\mathbf{P}[j_0 j_1 \cdots j_k]$ instead of $\mathbf{P}[\{\omega : j_0 j_1 \cdots j_k(\omega)\}]$. Assume that $j_{k-1} \in \{3,4,5,6,7,8\}$; we get

$$\mathbf{P}[j_0 j_1 \cdots j_k] = r_2 \mathbf{P}[j_0 j_1 \cdots j_{k-1}],$$

since $j_k$ is either 7 or 8 (depending on the value of $j_{k-1}$). If instead $j_{k-1} = 1$, then $j_k$ must be either $1, 3$ or $5$ and of course $j_l = 1$ for all $l \leq (k-1)$. Hence in this case

$$\mathbf{P}[j_0 j_1 \cdots j_k] = r_2 \mathbf{P}[j_0 j_1 \cdots j_{k-1}] = r_2 \mathbf{P}[\underbrace{11 \cdots 1}_{k}] = r_2 C r_1^{k+1}$$

if $j_k$ is equal to 3 or 5, and

$$\mathbf{P}[j_0 j_1 \cdots j_k] = \mathbf{P}[\underbrace{11 \cdots 1}_{k+1}] = C r_1^{k+2}$$

if $j_k = 1$. Similarly, if $j_{k-1} = 2$, then $j_k$ must be either $2, 4$ or $5$ and of course $j_l = 2$ for all $l \leq (k-1)$. Hence

$$\mathbf{P}[j_0 j_1 \cdots j_k] = r_2 \mathbf{P}[j_0 j_1 \cdots j_{k-1}] = r_2 \mathbf{P}[\underbrace{22 \cdots 2}_{k}] = r_2 C r_1^{k+1}$$

if $j_k$ is equal to 4 or 6, and

$$\mathbf{P}[j_0 j_1 \cdots j_k] = \mathbf{P}[\underbrace{22 \cdots 2}_{k+1}] = C r_1^{k+2}$$

if $j_k = 2$.

Let $\mathcal{A}_k$ be the set of all sequences $j_0 j_1 \cdots j_k$ corresponding to the event $h(Y_i, Y_{i+1}) = \cdots = h(Y_{i+k}, Y_{i+k+1}) = 1$. We have that

$$\mathbf{P}[h(Y_i, Y_{i+1}) = \cdots = h(Y_{i+k}, Y_{i+k+1}) = 1]$$

$$= \sum_{\mathcal{A}_k} \mathbf{P}[j_0 j_1 \cdots j_k]$$

$$= \sum_{\substack{\mathcal{A}_k \\ j_{k-1} \notin \{1,2\}}} \mathbf{P}[j_0 j_1 \cdots j_k] + \sum_{\substack{\mathcal{A}_k \\ j_{k-1} \in \{1,2\}}} \mathbf{P}[j_0 j_1 \cdots j_k]$$

$$= r_2 \sum_{\substack{\mathcal{A}_{k-1} \\ j_{k-1} \notin \{1,2\}}} \mathbf{P}[j_0 j_1 \cdots j_{k-1}] + 4 r_2 C r_1^{k+1} + 2 C r_1^{k+2}$$

$$= r_2 \left( \sum_{\substack{\mathcal{A}_{k-1} \\ j_{k-1} \notin \{1,2\}}} \mathbf{P}[j_0 j_1 \cdots j_{k-1}] + \mathbf{P}[\underbrace{11 \cdots 1}_{k}] + \mathbf{P}[\underbrace{22 \cdots 2}_{k}] \right)$$

$$+ 2 r_2 C r_1^{k+1} + 2 C r_1^{k+2}$$

$$= r_2 \sum_{\mathcal{A}_{k-1}} \mathbf{P}[j_0 j_1 \cdots j_{k-1}] + 2Cr_1^{k+1}(r_1 + r_2)$$

$$= r_2 \mathbf{P}[h(Y_i, Y_{i+1}) = \cdots = h(Y_{i+k-1}, Y_{i+k}) = 1] + r_1^{k+1}.$$

Comparing this to (9), and using $\mathbf{P}[h(Y_i, Y_{i+1}) = 1] = D_0$, we see that

$$\mathbf{P}[h(Y_i, Y_{i+1}) = \cdots = h(Y_{i+k}, Y_{i+k+1}) = 1] = D_k$$

for all $k \geq 0$, and so this case is also proved.

Finally the case $b > \frac{1}{2}$ remains. Take

$$g(x) = 1 - f(x) = 1 - b - 2|a|\cos(2\pi x - \phi) = 1 - b + 2|a|\cos(2\pi x - \phi'),$$

for some suitable choice of $\phi'$. Since $1 - b \leq \frac{1}{2}$, it follows from above that we can construct a two-block-factor $\{h(Y_i, Y_{i+1})\}_{i \in \mathbb{Z}}$ such that

$$\{h(Y_i, Y_{i+1})\}_{i \in \mathbb{Z}} \stackrel{\mathcal{D}}{=} \mathbf{P}^g.$$

With $\tilde{h} = 1 - h$, we get a new two-block-factor $\{\tilde{h}(Y_i, Y_{i+1})\}_{i \in \mathbb{Z}}$ with ones and zeros flipped. Lemma 2.4 in [9] then shows that $\{\tilde{h}(Y_i, Y_{i+1})\}_{i \in \mathbb{Z}}$ has distribution $\mathbf{P}^{1-g}$, which in turn is $\mathbf{P}^f$. □

When trying to generalize Theorem 1.3 to the case where $f$ is a trigonometric polynomial of degree $m$, one must consider not only the values of

$$\mathbf{P}^f[\eta(1) = \cdots = \eta(1+k) = 1],$$

but also the values of

$$\mathbf{P}^f[\eta(e_1) = \cdots = \eta(e_k) = 1],$$

where $e_i \in \mathbb{Z}\ \forall i \in \{1, \ldots, k\}$, but where $e_i$ is not necessarily equal to $e_{i-1} + 1$. Analyzing these new cylinder events adds to the complexity of the problem and therefore, in our opinion, the generalization of Theorem 1.3 (if indeed the generalization is true) does not seem to be trivial.

**Acknowledgment.** I would like to thank my supervisor Jeffrey Steif for presenting me the problem and for all the help I received during the writing of this paper.

DEPARTMENT OF MATHEMATICS
CHALMERS UNIVERSITY OF TECHNOLOGY
412 96 GOTHENBURG
SWEDEN
E-MAIL: broman@math.chalmers.se
URL: www.math.chalmers.se/~broman